\documentclass[12pt,a4paper]{article}
\usepackage{amsmath,amssymb,amsthm}

\def\N{\mathcal{N}}
\def\M{\mathcal{M}}

\def\add#1{{\operatorname{\textup{\textsf{add}}}}(#1)}

\def\cov#1{{\operatorname{\textup{\textsf{cov}}}}(#1)}

\def\fc{\mathfrak{c}}

\def\ssm{\smallsetminus}

\def\restrictedto%
{\mathchoice%
    {\!\upharpoonright\!}
    {\!\upharpoonright\!}
    {\upharpoonright}
    {\upharpoonright}
}

\def\dom{\operatorname{dom}}

\def\st{\mathchoice{:}{:}{\,:\,}{\,:\,}}

\def\size#1{\lvert#1\rvert}

\def\domn{\leq^*}

\newcommand{\splitnodes}[2][{}]%
	{\operatorname{\textup{\textsf{split}}}_{#1}({#2})}

\newcommand{\forces}[1][{}]{\Vdash_{#1}}
\newcommand{\forcestext}[2][{}]{\Vdash_{#1}\text{``}{#2}\text{''}}

\def\fl{\mathfrak{l}}

\theoremstyle{plain}
\newtheorem{thm}{Theorem}[section]
\newtheorem{lem}[thm]{Lemma}
\newtheorem{cor}[thm]{Corollary}

\theoremstyle{definition}
\newtheorem{defn}[thm]{Definition}
\theoremstyle{remark}

\newtheorem{remark}{Remark}

\begin{document}

\title{Covering a bounded set of functions 
	by an increasing chain of slaloms}
\author{Masaru Kada%
	\thanks{%
	Graduate School of Science, 
	Osaka Prefecture University. 
	Sakai, Osaka 599--8531 JAPAN}
	}
\date{}
\maketitle

\renewcommand{\thefootnote}{\relax}
\footnotetext[0]{2000 AMS Classification: Primary 03E17; Secondary 03E35.}
\renewcommand{\thefootnote}{\arabic{footnote}}

\begin{abstract}
A slalom is a sequence of finite sets of length $\omega$. 
Slaloms are ordered by coordinatewise inclusion 
with finitely many exceptions. 
Improving earlier results of Mildenberger, Shelah and Tsaban, 
we prove consistency results 
concerning 
existence and non-existence 
of an increasing sequence of a certain type of slaloms 
which covers a bounded set of functions in $\omega^\omega$. 
\end{abstract}

\section{Introduction}

We use standard terminology 
and refer the readers to \cite{BaJ:set} 
for undefined set-theoretic notions. 

Bartoszy\'nski \cite{Ba:comb} introduced 
the 
combinatorial concept 
of 
\emph{slalom} 
to study combinatorial aspects of 
measure and category on the real line. 

We call 
a sequence of 
finite subsets of $\omega$ of length $\omega$ 
a \emph{slalom}. 
For a function $g\in\omega^\omega$, 
let $\mathcal{S}^g$ be the set of slaloms $\varphi$ 
such that $\size{\varphi(n)}\leq g(n)$ for all $n<\omega$. 
$\mathcal{S}$ denotes $\mathcal{S}^g$ for $g(n)=2^n$. 
For 
two 
slaloms $\varphi$ and $\psi$, 
we write $\varphi\sqsubseteq\psi$ if 
$\varphi(n)\subseteq\psi(n)$ for all but finitely many $n<\omega$.
For a function $f\in\omega^\omega$ and a slalom $\varphi$, 
$f\sqsubseteq\varphi$ if 
$\langle \{f(n)\}\st n<\omega\rangle\sqsubseteq\varphi$.

Mildenberger, Shelah and Tsaban \cite{MiShTs:tau}
defined 
cardinals 
$\theta_h$ for $h\in\omega^\omega$ and $\theta_*$ 
to give a partial characterization of the cardinal $\mathfrak{od}$, 
the critical cardinality of a certain selection principle 
for 
open covers. 

The 
definition of $\theta_h$ 
in 
\cite{MiShTs:tau} 
is described using a combinatorial property 
which is called 
\emph{$o$-diagonalization}. 
Here we redefine $\theta_h$ 
to fit in the present context. 
It is easy to see that the following definition is equivalent to 
the original one. 
For a function 
$h\in(\omega\ssm\{0,1\})^\omega$, 
let $h-1$ denote the function $h'\in\omega^\omega$ 
which is defined by $h'(n)=h(n)-1$ for all $n$.

\begin{defn}
For a function 
$h\in(\omega\ssm\{0,1\})^\omega$, 
$\theta_h$ is the smallest size of a subset $\Phi$ of $\mathcal{S}^{h-1}$ 
which satisfies the following, \emph{if such a set\/ $\Phi$ exists}:
\begin{enumerate}
\item $\Phi$ is well-ordered by $\sqsubseteq$;
\item For every $f\in\prod_{n<\omega}h(n)$ there is a $\varphi\in\Phi$ 
	such that $f\sqsubseteq\varphi$. 
\end{enumerate}
If 
there is no such $\Phi$, 
we define $\theta_h=\fc^+$. 
\end{defn}

It is easy to see that $h_1\domn h_2$ implies $\theta_{h_1}\geq\theta_{h_2}$. 

\begin{defn}[\cite{MiShTs:tau}]
$\theta_*=\min\{\theta_h\st h\in\omega^\omega\}$. 
\end{defn}

In Section~\ref{sec:theta}, 
we will show that 
$\theta_*=\fc^+$ 
is consistent with ZFC. 

We say a proper forcing notion $\mathbb{P}$ 
has the \emph{Laver property} 
if, 
for any $h\in\omega^\omega$, 
$p\in\mathbb{P}$ and 
a $\mathbb{P}$-name $\dot{f}$ for a function in $\omega^\omega$ 
such that $p\forces[\mathbb{P}]{\dot{f}\in\prod_{n<\omega}h(n)}$, 
there exist $q\in\mathbb{P}$ and $\varphi\in\mathcal{S}$ 
such that $q$ is stronger than $p$ 
and $q\forces[\mathbb{P}]{\dot{f}\sqsubseteq\varphi}$. 

Mildenberger, Shelah and Tsaban
proved that $\theta_*=\aleph_1$ holds 
in 
all 
forcing 
models 
by 
a proper forcing notion with the Laver property 
over a model 
for CH, 
the continuum hypothesis 
\cite{MiShTs:tau}. 
In section~\ref{sec:theta}, 
we refine their result 
and state a sufficient condition 
for 
$\theta_*\leq\fc$. 
As a consequence, 
we will show that Martin's axiom implies 
$\theta_*=\fc$. 

In Section~\ref{sec:ht}, 
we give an application 
of 
the lemma presented in Section~\ref{sec:theta} 
to another problem in topology. 
We answer a question on 
approximations to 
the Stone--\v{C}ech compactification of $\omega$ 
by Higson compactifications 
of $\omega$, 
which was 
posed 
by Kada, Tomoyasu and Yoshinobu 
\cite{KTY:babylon}.

\section{Facts on the cardinal $\theta_*$}\label{sec:theta}

First 
we observe that 
$\theta_*=\fc^+$ 
is consistent with ZFC\@. 
We 
use the following theorem,  
which is 
a corollary 
of Kunen's classical result \cite{Ku:inacc}. 
For the readers' convenience, 
we 
present 
a 
complete 
proof in Section~\ref{sec:proof}.  

\begin{thm}\label{thm:isomorphismofnames}
Suppose that $\kappa\geq\aleph_2$. 
The following holds 
in the forcing model obtained by adding 
$\kappa$
Cohen reals 
over a model for 
CH:
Let $\mathcal{X}$ be a Polish space 
and $A\subseteq\mathcal{X}\times\mathcal{X}$ a 
Borel 
set. 
Then there is no sequence 
$\langle r_\alpha:\alpha<\omega_2\rangle$ 
in $\mathcal{X}$ 
which satisfies 
\[
\alpha\leq\beta<\omega_2 
\text{ if and only if } \langle r_\alpha,r_\beta\rangle\in A.	
\]
\end{thm}

%

Fix $h\in\omega^\omega$. 
We may regard $\mathcal{S}^{h-1}$ 
as a product space of countably many finite discrete spaces, 
and then the relation $\sqsubseteq$ 
on $\mathcal{S}^{h-1}$ 
is 
a Borel subset of $\mathcal{S}^{h-1}\times\mathcal{S}^{h-1}$. 

\begin{thm}\label{thm:thetainfty}
$\theta_*=\fc^+$ 
holds 
in the forcing model obtained by adding 
$\aleph_2$
Cohen reals 
over a model for 
CH. 
\end{thm}

\begin{proof}
Fix $h\in\omega^\omega$. 
By Theorem~\ref{thm:isomorphismofnames}, 
in the forcing model obtained by adding 
$\aleph_2$
Cohen reals 
over a model 
for 
CH, 
there is no $\sqsubseteq$-increasing chain of length $\omega_2$ 
in $\mathcal{S}^{h-1}$. 
This means that $\theta_h$ must be $\aleph_1$ 
whenever $\theta_h\leq\fc$. 

On the other hand, 
$\cov\M=\aleph_2$ holds in the same model. 
Also, 
by 
\cite{MiShTs:tau} 
we have $\cov\M\leq\mathfrak{od}\leq\theta_h$.  
This means that $\theta_h$ cannot be $\aleph_1$ in this model, 
and hence 
$\theta_h=\fc^+$. 
\end{proof}

Next we state a sufficient condition for 
$\theta_*\leq\fc$. 
%
We 
use 
the following 
characterization 
of $\add\N$. 

\begin{thm}[{\cite[Theorem~2.3.9]{BaJ:set}}]\label{thm:addn}
$\add\N$ is the smallest size of a subset $F$ of $\omega^\omega$ 
such that, 
for every $\varphi\in\mathcal{S}$ 
there is an $f\in F$ such that $f\not\sqsubseteq\varphi$. 
\end{thm}


\begin{defn}[{\cite[Section~5]{Kada:gamecd}}]
For a function $h\in\omega^\omega$, 
$\mathfrak{l}_h$ 
is the smallest size of a subset $\Phi$ of $\mathcal{S}$ 
such that 
for all $f\in\prod_{n<\omega}h(n)$ 
there is a $\varphi\in\Phi$ 
such that $f\sqsubseteq\varphi$. 
Let $\mathfrak{l}=\sup\{\mathfrak{l}_h\st h\in\omega^\omega\}$. 
\end{defn}

Note that 
$h_1\leq^* h_2$ implies $\mathfrak{l}_{h_1}\leq\mathfrak{l}_{h_2}$. 


If 
CH 
holds in a ground model $V$, 
$h\in\omega^\omega\cap V$, 
and a proper forcing notion $\mathbb{P}$ has the Laver property, 
then $\mathfrak{l}_h=\aleph_1$ holds in the model $V^{\mathbb{P}}$. 
Consequently, 
if 
CH 
holds in $V$, 
$\langle\mathbb{P}_\alpha,\dot{\mathbb{Q}}_\alpha\st\alpha<\omega_2\rangle$ 
is a countable support iteration 
of proper forcings, 
$\mathbb{P}=\lim_{\alpha<\omega_2}\mathbb{P}_\alpha$  
and 
\[
\forcestext[{\mathbb{P}_\alpha}]{%
	\size{\dot{\mathbb{Q}}_\alpha}\leq\aleph_1
	\text{ and }
	\dot{\mathbb{Q}}_\alpha\text{ has the Laver property}
	}
\] 
holds for every $\alpha<\omega_2$, 
then $\mathfrak{l}=\aleph_1$ holds in $V^{\mathbb{P}}$, 
since 
every function $h$ in $V^\mathbb{P}$ 
appears in $V^{\mathbb{P}_\alpha}$ for some $\alpha<\omega_2$, 
where CH holds.\footnote{%
In the paper \cite{KTY:babylon}, 
the authors 
state 
``If CH holds in a ground model $V$, 
and a proper forcing notion $\mathbb{P}$ has the Laver property, 
then $\mathfrak{l}=\aleph_1$ holds in the model $V^{\mathbb{P}}$''. 
But it is inaccurate, 
since we do not see 
the values of $\mathfrak{l}_h$ 
for functions $h\in V^{\mathbb{P}}$ 
which are not bounded by any function from $V$. 
}

Now we define a subset $\mathcal{S}^+$ of $\mathcal{S}$ 
as follows: 
\[
\mathcal{S}^+=\left\{\varphi\in\mathcal{S}\st
	\lim_{n\to\infty}
	\frac{\size{\varphi(n)}}{2^n}
	=0\right\}. 
\]
Let $\fl'_{h}$ 
be the smallest size of a subset $\Phi$ of $\mathcal{S}^+$ 
such that 
for all $f\in\prod_{n<\omega}h(n)$ 
there is a $\varphi\in\Phi$ 
such that $f\sqsubseteq\varphi$. 
Clearly we have $\fl_h\leq\fl'_h$, 
and 
it is easy to see that 
for every $h\in\omega^\omega$ 
there is an $h^*\in\omega^\omega$ 
such that $\fl'_h\leq\fl_{h^*}$. 
Hence we have $\fl=\sup\{\fl'_h\st h\in\omega^\omega\}$.

\begin{lem}
For a subset $\Phi$ of $\mathcal{S}^+$ of size less than 
$\add\N$, 
there is a $\psi\in\mathcal{S}^+$ 
such that 
$\varphi\sqsubseteq\psi$ for all $\varphi\in\Phi$. 
\end{lem}

\begin{proof}
For each $\varphi\in\mathcal{S}^+$, 
define 
an increasing function 
$\eta_\varphi\in\omega^\omega$ 
by letting 
\[
	\eta_\varphi(m)=\min\left\{l<\omega\st
		\forall k\geq l\,\,
		\left(\,\size{\varphi(k)}<\frac{2^k}{m\cdot2^m}\,\right)
	\right\}
\]
for all $m<\omega$. 
$\eta_\varphi$ is well-defined by the definition of $\mathcal{S}^+$. 

Suppose $\kappa<\add\N$ 
and fix a set $\Phi\subseteq\mathcal{S}^+$ of size $\kappa$ 
arbitrarily. 
Since $\kappa<\add\N\leq\mathfrak{b}$, 
there is 
a 
function $\eta\in\omega^\omega$ 
such that 
$\lim_{n\to\infty}{\eta(n)}/{2^n}=\infty$  
and 
for all $\varphi\in\Phi$ 
we have $\eta_\varphi\leq^*\eta$. 
For each $m<\omega$, 
let 
$I_m=\{\eta(m),\eta(m)+1,\ldots,\eta(m+1)-1\}$ 
and enumerate 
$
\prod_{n\in I_m}[\omega]^{\leq\lfloor 2^n/(m\cdot 2^m)\rfloor}
$
as $\{s_{m,i}: i<\omega\}$,	
where $\lfloor r\rfloor$ denotes 
the largest integer which does not exceed 
the real number $r$. 

For $\varphi\in\Phi$, 
define 
$\tilde{\varphi}\in\omega^\omega$ 
as follows. 
If there is an $i<\omega$ 
such that $\varphi\mathbin{\upharpoonright}I_m=s_{m,i}$, 
then let $\tilde{\varphi}(m)=i$; 
otherwise $\tilde{\varphi}(m)$ is arbitrary. 

Since $\size{\Phi}=\kappa<\add\N$ 
and by Theorem~\ref{thm:addn}, 
there is a 
$\hat{\psi}\in\mathcal{S}$ 
such that, 
for all $\varphi\in\Phi$ 
we have 
$\tilde{\varphi}\sqsubseteq\hat{\psi}$. 
Define $\psi$ 
by letting for each $n$, 
if $n\in I_m$ then 
$\psi(n)=\bigcup\{s_{m,i}(n):i\in\hat{\psi}(m)\}$, 
and 
if $n<\eta(0)$ then $\psi(n)=\emptyset$. 
It is straightforward 
to check that 
$\psi\in\mathcal{S}^+$ 
and $\varphi\sqsubseteq\psi$ for all $\varphi\in\Phi$. 
\end{proof}

\begin{lem}\label{lem:addnl}
Suppose that 
$h\in\omega^\omega$ 
satisfies 
$h(n)>n^2$ for all $n<\omega$. 
If 
$\add\N=\fl'_h=\kappa$, 
then 
there is 
an $\sqsubseteq$-increasing sequence 
$\langle\sigma_\alpha\st\alpha<\kappa\rangle$ 
in 
$\mathcal{S}^+$
such that, 
for all $f\in\prod_{n<\omega}h(n)$ there is an $\alpha<\kappa$ 
such that $f\sqsubseteq\varphi_\alpha$. 
\end{lem}

\begin{proof}
Fix a sequence 
$\langle \varphi_\alpha:\alpha<\kappa\rangle$ 
in 
$\mathcal{S}^+$ 
so that 
for all $f\in\prod_{n<\omega}h(n)$ 
there is an $\alpha<\kappa$ 
such that $f\sqsubseteq\varphi_\alpha$. 
Using the previous lemma, 
inductively construct an $\sqsubseteq$-increasing sequence 
$\langle\sigma_\alpha:\alpha<\kappa\rangle$ 
of elements of $\mathcal{S}^+$ 
so that 
$\varphi_\alpha\sqsubseteq\sigma_\alpha$ 
holds 
for each $\alpha<\omega_2$. 
Then 
$\langle\sigma_\alpha:\alpha<\kappa\rangle$ 
is as required. 
\end{proof}

Define $H_1\in\omega^\omega$ by letting $H_1(n)=2^n+1$ for all $n$.

\begin{thm}\label{thm:addnltheta}
If 
$\add\N=\mathfrak{l}'_{H_1}$, 
then $\theta_*=\mathfrak{od}=\add\N$. 
\end{thm}

\begin{proof}
Let $\kappa=\add\N=\mathfrak{l}'_{H_1}$. 
Since $\mathcal{S}^+\subseteq\mathcal{S}\subseteq\mathcal{S}^{H_1-1}$, 
the previous lemma 
shows that 
$\theta_*\leq\theta_{H_1}\leq\kappa$. 
On the other hand, 
by \cite{MiShTs:tau}, 
we have 
$\kappa=\operatorname{\textsf{add}}(\mathcal{N})
	\leq\operatorname{\textsf{cov}}(\mathcal{M})
	\leq\mathfrak{od}
	\leq\theta_*$. 
\end{proof}

\begin{cor}[\cite{MiShTs:tau}]
If a ground model\/ $V$ satisfies 
CH, 
and a proper forcing notion $\mathbb{P}$ 
has the Laver property, 
then $\theta_*=\aleph_1$ holds in the model\/ $V^{\mathbb{P}}$. 
\end{cor}

\begin{proof}
Follows from Theorem~\ref{thm:addnltheta} and the fact that 
$\operatorname{\textsf{add}}(\mathcal{N})
=\mathfrak{l}'_{H_1}=\mathfrak{l}_{H_1^*}=\aleph_1$ holds 
in the model $V^{\mathbb{P}}$. 
\end{proof}

\begin{cor}
Martin's axiom implies $\theta_*=\fc$. 
\end{cor}

\begin{proof}
Follows from Theorem~\ref{thm:addnltheta} and the fact that 
$\add\N=\mathfrak{l}'_{H_1}=\mathfrak{l}=\fc$ holds 
under Martin's axiom. 
\end{proof}

\section{Application
	}\label{sec:ht}

In this section, 
we give an answer to a question 
which was 
posed 
by Kada, Tomoyasu and Yoshinobu \cite{KTY:babylon}. 
We refer the reader to $\cite{KTY:babylon}$ 
for undefined topological notions. 

For compactifications $\alpha X$ and $\gamma X$ 
of a completely regular Hausdorff space $X$, 
we write $\alpha X\leq\gamma X$ 
if there is a continuous surjection from $\gamma X$ to $\alpha X$ 
which fixes the points from $X$, 
and $\alpha X\simeq\gamma X$ 
if $\alpha X\leq\gamma X\leq\alpha X$. 
The Stone--\v{C}ech compactification $\beta X$ 
of $X$
is the maximal compactification of $X$ 
in the sense of the order relation $\leq$ among compactifications of $X$. 


For a proper metric space $(X,d)$, 
$\overline{X}^d$ denotes the Higson compactification 
of $X$ with respect to the metric $d$. 

$\mathfrak{ht}$ is the smallest size of a set $D$ 
of proper metrics on $\omega$ such that
\begin{enumerate}
\item $\{\overline{\omega}^d\st d\in D\}$ is well-ordered by $\leq$; 
\item There is no $d\in D$ such that 
	$\overline{\omega}^d\simeq\beta\omega$; 
\item $\beta\omega\simeq\sup\{\overline{\omega}^d\st d\in D\}$, 
	where $\sup$ is in the sense of the order relation $\leq$ 
	among compactifications of $\omega$; 
\end{enumerate}
\emph{if such a set $D$ exists}. 
We define $\mathfrak{ht}=\fc^+$ 
if there is no such $D$. 

Kada, Tomoyasu and Yoshinobu 
\cite[Theorem~6.16]{KTY:babylon}
proved 
the consistency of 
$\mathfrak{ht}=\fc^+$ 
using a similar argument to the proof of Theorem~\ref{thm:thetainfty}. 
But the consistency 
of 
$\mathfrak{ht}\leq\fc$ 
was not addressed. 
Here 
we state a sufficient condition for 
$\mathfrak{ht}\leq\fc$, 
and show that it is consistent with ZFC\@. 

Define $H_2\in\omega^\omega$ by letting $H_2(n)=2^{2^{(n^4)}}$ for all $n$. 
The following lemma 
is obtained 
as a corollary of the proof of \cite[Theorem~6.11]{KTY:babylon}. 

\begin{lem}\label{lem:scsl}
Let $\kappa$ be a cardinal. 
If there is 
an $\sqsubseteq$-increasing sequence 
$\langle\varphi_\alpha\st\alpha<\kappa\rangle$ 
of slaloms in $\mathcal{S}$ such that 
for all $f\in\prod_{n<\omega}H_2(n)$ there is an $\alpha<\kappa$ 
such that $f\sqsubseteq\varphi_\alpha$, 
then $\mathfrak{ht}\leq\kappa$. 
\end{lem}

Now we have the following theorem. 

\begin{thm}\label{thm:addnlht}
If $\add\N=\fl'_{H_2}$, then $\mathfrak{ht}=\add\N$. 
\end{thm}

\begin{proof}
$\add\N\leq\mathfrak{ht}$ 
is proved in \cite[Section~6]{KTY:babylon}. 
To see 
$\mathfrak{ht}\leq\add\N$, 
apply Lemma~\ref{lem:addnl} 
for $h=H_2$ 
to get a sequence of slaloms 
which is required in Lemma~\ref{lem:scsl}. 
\end{proof}

\begin{cor}
If a ground model\/ $V$ satisfies 
CH, 
and a proper forcing notion $\mathbb{P}$ 
has the Laver property, 
then $\mathfrak{ht}=\aleph_1$ holds in the model\/ $V^{\mathbb{P}}$. 
\end{cor}

\begin{proof}
Follows from Theorem~\ref{thm:addnlht} and the fact that 
$\operatorname{\textsf{add}}(\mathcal{N})
	=\mathfrak{l}'_{H_2}=\mathfrak{l}_{H_2^*}=\aleph_1$ holds 
in the model $V^{\mathbb{P}}$. 
\end{proof}

\begin{cor}
Martin's axiom implies $\mathfrak{ht}=\fc$. 
\end{cor}

\begin{proof}
Follows from Theorem~\ref{thm:addnlht} and the fact that 
$\add\N=\mathfrak{l}'_{H_2}=\mathfrak{l}=\fc$ holds 
under Martin's axiom. 
\end{proof}

\section{Proof of Theorem~\ref{thm:isomorphismofnames}}
\label{sec:proof}

This section 
is devoted 
to the proof of Theorem~\ref{thm:isomorphismofnames}. 
The idea of the proof 
is the same as the one in Kunen's 
original proof \cite{Ku:inacc}, 
which is known as the ``isomorphism of names'' argument. 
The same argument is also found in \cite{JSS:combprinciple}. 

For an infinite set $I$, 
let $\mathbb{C}(I)=\operatorname{Fn}(I,2,\aleph_0)$, 
the 
canonical 
Cohen forcing notion 
for 
the index set $I$. 
As described in \cite[Chapter~7]{Ku:set}, 
for any $\mathbb{C}(I)$-name $\dot{r}$ for a subset of $\omega$, 
we can find a countable subset $J$ of $I$ 
and a \emph{nice $\mathbb{C}(J)$-name} $\dot{s}$ for a subset of $\omega$ 
such that $\forces[\mathbb{C}(I)]{\dot{s}=\dot{r}}$. 
For a countable set $I$, 
there are only 
$\mathfrak{c}$ 
nice $\mathbb{C}(I)$-names 
for subsets of $\omega$. 




\begin{proof}[Proof of Theorem~\ref{thm:isomorphismofnames}]
Suppose that $\kappa\geq\aleph_2$. 
Let $\mathcal{X}$ be a Polish space, 
$\dot{A}$ a $\mathbb{C}(\kappa)$-name 
for a Borel subset of $\mathcal{X}\times\mathcal{X}$, 
%
%
and 
$\langle\dot{r}_\alpha\st\alpha<\omega_2\rangle$ 
a sequence 
of $\mathbb{C}(\kappa)$-names for elements of $\mathcal{X}$. 

We will prove the following statement:
\[
	\forces[\mathbb{C}(\kappa)]{
	\exists\alpha<\omega_2\,
	\exists\beta<\omega_2\,
	(\alpha<\beta\land(
	\langle\dot{r}_\alpha,\dot{r}_\beta\rangle\notin\dot{A}
	\lor
	\langle\dot{r}_\beta,\dot{r}_\alpha\rangle\in\dot{A})
	).
	}
\]
There is nothing to do if it holds that 
\[
	\forces[\mathbb{C}(\kappa)]{
	\exists\alpha<\omega_2\,
	\exists\beta<\omega_2\,
	(\alpha<\beta\land
	\langle\dot{r}_\alpha,\dot{r}_\beta\rangle\notin\dot{A}).
	}
\]
So we assume that it fails, 
and fix any $p\in\mathbb{C}(\kappa)$ which satisfies 
\begin{equation}
	p
	\forces[\mathbb{C}(\kappa)]{	
	\forall\alpha<\omega_2\,
	\forall\beta<\omega_2\,
	(\alpha<\beta\rightarrow
	\langle\dot{r}_\alpha,\dot{r}_\beta\rangle\in\dot{A}).}
	\tag*{$(*)$}
\end{equation}
We 
will 
find 
$\alpha,\beta<\omega_2$ 
such that 
$\alpha<\beta$ and 
$
	p
	\forces[\mathbb{C}(\kappa)]{
	\langle\dot{r}_\beta,\dot{r}_\alpha\rangle\in\dot{A}}
$, 
which concludes the proof. 

Let $J_p=\dom(p)$. 
Find a set $J_A\in[\kappa]^{\aleph_0}$ 
and a nice $\mathbb{C}(J_A)$-name $\dot{C}_A$ 
for a subset of $\omega$ 
such that 
\[
	\forcestext[\mathbb{C}(\kappa)]{
		\dot{C}_A\text{ is a Borel code of }\dot{A}.} 
\]
For each $\alpha<\omega_2$, 
find 
a set 
$J_\alpha\in[\kappa]^{\aleph_0}$ 
and a nice 
$\mathbb{C}(J_\alpha)$-name $\dot{C}_\alpha$ 
for a subset of $\omega$ 
such that 
\[
	\forcestext[\mathbb{C}(\kappa)]{
		\dot{C}_\alpha
		\text{ is a Borel code of }
		\{\dot{r}_\alpha\}.
	}
\]
Using the $\Delta$-system lemma \cite[II Theorem~1.6]{Ku:set}, 
take $S\in[\kappa]^{\aleph_0}$ 
and $K\in[\omega_2]^{\aleph_2}$ 
so that 
$
J_p\cup 
J_A\cup(J_\alpha\cap J_\beta)\subseteq S
$ 
for any $\alpha,\beta\in K$ with $\alpha\neq\beta$.  
Without loss of generality 
we may assume that 
$\size{J_\alpha\ssm S}=\aleph_0$ for all $\alpha\in K$. 
For each $\alpha\in K$, 
enumerate $J_\alpha\ssm S$ as 
$\langle\delta^\alpha_n\st n<\omega\rangle$. 

For $\alpha,\beta\in K$, 
and let $\sigma_{\alpha,\beta}$ 
be the involution (automorphism of order 2) 
of $\mathbb{C}(\kappa)$ 
obtained by the permutation of coordinates 
which interchanges $\delta^\alpha_n$ with $\delta^\beta_n$ 
for each $n$. 
$\sigma_{\alpha,\beta}$ naturally induces 
an 
involution 
of the class of all $\mathbb{C}(\kappa)$-names: 
We simply denote it by $\sigma_{\alpha,\beta}$. 
%
Since 
$
J_p\cup 
J_A\subseteq S
$, 
for all $\alpha,\beta\in K$ 
we have 
$\sigma_{\alpha,\beta}(p)=p$, 
$\sigma_{\alpha,\beta}(\dot{C}_A)=\dot{C}_A$ 
and 
$
	\forces[\mathbb{C}(\kappa)]{
		\sigma_{\alpha,\beta}(\dot{A})=\dot{A}}
$. 



Since $\size{K}=\aleph_2$ 
and there are only 
$\mathfrak{c}=\aleph_1$ 
nice names for subsets of $\omega$ over a countable index set, 
we can find $\alpha,\beta\in K$ 
with 
$\alpha<\beta$ 
such that  
$\sigma_{\alpha,\beta}(\dot{C}_\alpha)=\dot{C}_\beta$. 
Then 
$\sigma_{\alpha,\beta}(\dot{C}_\beta)
=\dot{C}_\alpha$ 
and 
\[
	\forcestext[\mathbb{C}(\kappa)]{
	\sigma_{\alpha,\beta}(\dot{r}_\alpha)=\dot{r}_\beta
	\text{ and }
	\sigma_{\alpha,\beta}(\dot{r}_\beta)=\dot{r}_\alpha. }
\] 

By $(*)$,  
we have 
$p\forces[\mathbb{C}(\kappa)]
	{\langle\dot{r}_\alpha,\dot{r}_\beta\rangle}\in\dot{A}$. 
Since $\sigma_{\alpha,\beta}$ is an automorphism of $\mathbb{C}(\kappa)$, 
we have 
\[
	\sigma_{\alpha,\beta}(p)\forces[\mathbb{C}(\kappa)]{
	\langle\sigma_{\alpha,\beta}(\dot{r}_\alpha),
		\sigma_{\alpha,\beta}(\dot{r}_\beta)\rangle
	\in\sigma_{\alpha,\beta}(\dot{A})}	
\]
and hence 
$p\forces[\mathbb{C}(\kappa)]
	{\langle\dot{r}_\beta,\dot{r}_\alpha\rangle}\in\dot{A}$. 
\end{proof}

\begin{remark}
Fuchino pointed out that 
Theorem~\ref{thm:isomorphismofnames} 
is generalized 
in the following 
two ways \cite{BrFu:colorord}: 
(1)~ 
The set $A$ is not necessarily Borel, but is 
``definable'' by some formula. 
(2)~
We can prove a similar result 
for a forcing extension by a side-by-side product 
of the same forcing notions, 
each generically adds a real in a natural way. 
The argument in the above proof 
also works in those generalized settings.
\end{remark}

\subsection*{Acknowledgement}

I thank 
Saka\'e Fuchino, 
Hiroshi Fujita, 
Teruyuki Yorioka 
and the referee 
for 
helpful 
comments and remarks 
on Theorem~2.1 and its proof.


\begin{thebibliography}{1}

\bibitem{Ba:comb}
T.~Bartoszy\'nski.
\newblock Combinatorial aspects of measure and category.
\newblock {\em Fund. Math.}, 127:225--239, 1987.

\bibitem{BaJ:set}
T.~Bartoszy\'nski and H.~Judah.
\newblock {\em Set Theory: On the Structure of the Real Line}.
\newblock A. K. Peters, Wellesley, Massachusetts, 1995.

\bibitem{BrFu:colorord}
J.~Brendle and S.~Fuchino. 
\newblock Coloring ordinals by reals. 
\newblock preprint. 


\bibitem{JSS:combprinciple}
I.~Juh\'asz, L.~Soukup, and Z.~Szentmikl\'ossy.
\newblock Combinatorial principles from adding {Cohen} reals.
\newblock In J.~A. Makowsky, editor, {\em Logic Colloquium 95, Proceedings of
  the Annual European Summer Meeting of the Association of Symbolic Logic},
  Lecture Notes in Logic. 11., pages 79--103, Haifa, Israel, 1998. Springer.
\bibitem{Kada:gamecd}
M.~Kada.
\newblock More on {Cicho\'n's} diagram and infinite games.
\newblock {\em J. Symbolic Logic}, 65:1713--1724, 2000.

\bibitem{KTY:babylon}
M.~Kada, K.~Tomoyasu, and Y.~Yoshinobu.
\newblock How many miles to $\beta\omega$? --- {Approximating} $\beta\omega$ by
  metric-dependent compactifications.
\newblock {\em Topology Appl.}, 145:277--292, 2004.

\bibitem{Ku:inacc}
K.~Kunen.
\newblock {\em Inaccessibility properties of cardinals}.
\newblock {Ph.D.} dissertation, Stanford, 1968.

\bibitem{Ku:set}
K.~Kunen.
\newblock {\em Set Theory: an introduction to independence proofs}, volume 102
  of {\em Studies in Logic}.
\newblock North Holland, 1980.

\bibitem{MiShTs:tau}
H.~Mildenberger, S.~Shelah, and B.~Tsaban.
\newblock The combinatorics of $\tau$-covers.
\newblock {\em Topology Appl.}, to appear.

\end{thebibliography}

\end{document}